\documentclass[12pt]{article}

\usepackage{amsthm, amssymb, amstext}
\usepackage{comment}
\usepackage{setspace}
\usepackage[margin=2.5cm]{geometry}
\newtheorem{theorem}{Theorem}

\theoremstyle{definition}

\newtheorem{remark*}{Remark}

\usepackage{amsthm}

\title{Paths between colourings of sparse graphs}
\author{Carl Feghali\thanks{Department of Informatics, University of Bergen, Norway, email: \texttt{feghali.carl@gmail.com} }}
\date{}

\begin{document}
\maketitle

\begin{abstract}
The reconfiguration graph $R_k(G)$ of the $k$-colourings
of a graph~$G$ has as vertex set the set of all possible $k$-colourings
of $G$ and two colourings are adjacent if they differ on the colour of exactly
one vertex.  

We give a short proof of the following theorem of Bousquet and Perarnau (\emph{European Journal of Combinatorics}, 2016). Let $d$ and $k$ be positive integers, $k \geq d+1$. For every $\epsilon > 0$ and every graph $G$ with $n$ vertices and maximum average degree $d - \epsilon$, there exists a constant $c = c(d, \epsilon)$ such that $R_k(G)$ has diameter $O(n^c)$.  Our proof can be transformed into a simple polynomial time algorithm that finds a path between a given pair of colourings in $R_k(G)$. 
\end{abstract}

Let $G$ be a graph, and let $k$ be a non-negative integer. 
A $k$-colouring of $G$ is a function $f: V(G) \rightarrow \{1, \dots, k\}$ such that $f(u) \not= f(v)$ whenever $(u, v) \in E(G)$. 
The reconfiguration graph $R_k(G)$ of the $k$-colourings of $G$ has as vertex set the set of all $k$-colourings of $G$ and two vertices of $R_k(G)$ are adjacent if they differ on the colour of exactly one vertex. 

The first results on reconfiguration graphs appeared in the context of the irreducibility of Glauber dynamics for graph colourings; see \cite{jerrum, vigoda}. It was later studied by Cereceda, van den Heuvel and Johnson \cite{CHJ06, CHJ06a} who proved some results on the connectivity of reconfiguration graphs. Bonsma and Cereceda \cite{BC09} also considered bounding the diameter of (the components of) $R_k(G)$ for arbitrary graphs $G$ by obtaining a polynomial bound (in the order of $G$) whenever $k = 3$ and a superpolynomial bound (in the size of $G$) whenever $k \geq 4$. This line of research was continued by other researchers where a polynomial bound on the diameter of $R_k(G)$ is sought for every $k \geq 4$ and every graph $G$ belonging to some special graph class such as, for example, chordal and chordal bipartite graphs \cite{BJLPP14}, graphs of bounded treewidth  \cite{bonamy13} and graphs of bounded maximum average degree \cite{bousquet11}. Reconfiguration problems corresponding to decision problems other than graph colourings have also been the subject of much attention; see \cite{nishimura, He13} for excellent surveys on the topic. 

Our purpose in this note is to give a significantly shorter proof of the aforementioned result in \cite{bousquet11} due to Bousquet and Perarnau. 

\begin{theorem}\label{mainthm}
Let $d$ and $k$ be positive integers, $k \geq d+1$. For every $\epsilon > 0$ and every graph $G$ with $n$ vertices and maximum average degree $d - \epsilon$, there exists a constant $c = c(d, \epsilon)$ such that $R_k(G)$ has diameter $O(n^c)$. 
\end{theorem} 

Let us emphasise that our most important message is not the length of the proof but the technique used to establish Theorem \ref{mainthm}. Roughly speaking, previous approaches to prove results of this kind involve transforming, one vertex at a time, initial configurations to `more structured' intermediate configurations, as well as reducing the order of the graph by a small fraction of its order. We shall take the reverse route of directly transforming one configuration into the other, simultaneously reconfiguring several vertices at a time (even with the worst possible sequence of reconfiguration steps) while ensuring that the order of the graph can be reduced by a large fraction of its order. This way the overall number of reconfiguration steps can still be bounded by a polynomial function of the order of the graph. 

Our proof also has the benefit of being straightforwardly adaptable into an algorithm that is much simpler than the one in \cite{bousquet11}. That is, given $\epsilon > 0$, $k\geq d+1$, a graph $G$ with maximum average degree $d - \epsilon$ and $k$-colourings $\alpha$ and $\beta$ of $G$, one can find a path between $\alpha$ and $\beta$ in $R_k(G)$ in time that is polynomial in the order of $G$.

\begin{proof}[Proof of Theorem \ref{mainthm}] Recall that the maximum average degree of $G$ is defined as
\[
\max \bigg\{ \frac{2|E(H)|}{|V(H)|} : H \subseteq G \bigg\}. 
\]

Let $H$ be any subgraph of $G$, and let $h = |V(H)|$. An independent set $I$ of $H$ is said to be \emph{special} if it contains at least $\epsilon h / d^2$ vertices and every vertex of $I$ has at most $d - 1$ neighbours in $H$. Let us first show that $H$ contains  a special independent set. 

Let $S$ be the set of vertices of degree $d - 1$ or less in $H$. The number of vertices of $S$ is at least $\epsilon h/ d$ since otherwise
\[
\sum_{v \in H} \textrm{deg}(v) \geq \sum_{v \in H - S} \textrm{deg}(v) > d\bigg(h - \frac{\epsilon h}{d}\bigg) = (d - \epsilon)h,
\]
which contradicts the maximum average degree of $G$. Let $I \subseteq S$ be a maximal independent set of $S$. Then every vertex of $S - I$ has at least one neighbour in $I$ and every vertex of $I$ has at most $d - 1$ neighbours in $S$. Therefore, $|I| + (d- 1)|I| \geq |S| \geq \epsilon h / d$ and so $I$ is a special independent set, as required. 

In order to complete the proof of the theorem, it suffices to show that we can recolour from any $k$-colouring $\alpha$ of $G$ to any other $k$-colouring $\beta$ by $O(n^c)$ recolourings for some constant $c = c(d, \epsilon)$. Let $I$ be a special independent of $G$, and let $H = G - I$. Let $\alpha^H$ and $\beta^H$ denote, respectively, the restrictions of $\alpha$ and $\beta$ to $H$. Assume that there is a sequence of recolourings from $\alpha^H$ to $\beta^H$. 

We can extend this sequence to $H + u$ for each $u \in I$ by recolouring $u$ with some colour not used in its neighbourhood (this is possible since $k \geq d+1$) whenever some neighbour of $u$ is recoloured to its colour. At the end of the sequence, we recolour $u$ with $\beta(u)$. So we can successively add vertices of $I$ to $H$ and obtain in this way a sequence of recolourings from $\alpha$ to $\beta$ such that the number of recolourings of vertices of $I$ are independent of one another. As every subgraph of $G$ contains a special independent set, such a sequence exists (say, by strong induction on the order of $G$). Moreover, the maximum number $f(n)$ of times  each vertex of a graph with $n$ vertices and maximum average degree $d - \epsilon$ is recoloured satisfies the recurrence relation
\[
f(n) \leq (d - 1) \cdot f\bigg(n - \frac{\epsilon}{d^2}n\bigg) + 1,
\]
and the theorem follows by the master theorem. 
\end{proof}

\begin{remark*}
Although our proof is shorter than the one in \cite{bousquet11}, our bound on the the diameter of the reconfiguration graph is worse. More specifically, in our case $c = \log_b(a) + 1$ where $a = d- 1$ and $b = d^2/(d^2 - \epsilon)$ whereas  in \cite{bousquet11} one has $c= \frac{1}{ \log_d (d/(d - \epsilon))} + 2$. 
\end{remark*}

\section*{Acknowledgements}

The author is grateful to the referees for carefully checking the paper and to a referee whose comments have helped to improve the exposition of the paper. This work was supported by the research Council of Norway via the project CLASSIS, grant number 249994.

 \bibliography{bibliography}{}
\bibliographystyle{abbrv}
 
\end{document}